\documentclass[times,preprint]{elsarticle}

\usepackage{lineno}
 \RequirePackage{geometry}
 \global\bibsep=0pt

\usepackage{framed,multirow}

\journal{NumHyp2019 Proceedings}
\usepackage{type1cm}        
\usepackage{makeidx}         
\usepackage{graphicx}        
\usepackage{multicol}        
\usepackage[bottom]{footmisc}

\usepackage{newtxtext}       %
\usepackage{newtxmath}       

\makeindex             

\usepackage{lipsum}
\usepackage{amsfonts}
\usepackage{graphicx}
\usepackage{epstopdf}
\usepackage{algorithmic}
\ifpdf
\DeclareGraphicsExtensions{.eps,.pdf,.png,.jpg}
\else
\DeclareGraphicsExtensions{.eps}
\fi

\usepackage{multicol}        
\usepackage[bottom]{footmisc}

\usepackage{amssymb}
\usepackage{amsmath}
\usepackage{wasysym}

\usepackage{bm}
\usepackage{xspace}

\usepackage{booktabs}
\usepackage{tabularx}

\usepackage[labelfont=bf]{caption}
\usepackage{subfig}

\usepackage{color}

\usepackage{hyperref}
\usepackage[capitalize]{cleveref}

\crefname{section}{Sect.}{Sect.}

\usepackage{xargs}  

\renewcommand{\vec}[1]{{\bf #1}}

\def\@thmcountersep{.}

\newcolumntype{L}[1]{>{\raggedright\arraybackslash}p{#1}} 
\newcolumntype{C}[1]{>{\centering\arraybackslash}p{#1}} 
\newcolumntype{R}[1]{>{\raggedleft\arraybackslash}p{#1}} 
\newcolumntype{Y}{>{\centering\arraybackslash}X} 
\newcolumntype{Z}{>{\raggedleft\arraybackslash}X} 

\def\@thmcountersep{.}

\makeatletter
\newcommand{\manuallabel}[2]{\def\@currentlabel{#2}\label{#1}}
\makeatother

\newcounter{subcount}

\DeclareMathAlphabet{\mathitbf}{OML}{cmm}{b}{it}

\makeatletter
\newcommand\RedeclareMathOperator{%
	\@ifstar{\def\rmo@s{m}\rmo@redeclare}{\def\rmo@s{o}\rmo@redeclare}%
}
\newcommand\rmo@redeclare[2]{%
	\begingroup \escapechar\m@ne\xdef\@gtempa{{\string#1}}\endgroup
	\expandafter\@ifundefined\@gtempa
	{\@latex@error{\noexpand#1undefined}\@ehc}%
	\relax
	\expandafter\rmo@declmathop\rmo@s{#1}{#2}}
\newcommand\rmo@declmathop[3]{%
	\DeclareRobustCommand{#2}{\qopname\newmcodes@#1{#3}}%
}
\@onlypreamble\RedeclareMathOperator
\makeatother

\DeclareMathSymbol{\alpha}{\mathalpha}{letters}{"0B}
\DeclareMathSymbol{\beta}{\mathalpha}{letters}{"0C}
\DeclareMathSymbol{\gamma}{\mathalpha}{letters}{"0D}
\DeclareMathSymbol{\delta}{\mathalpha}{letters}{"0E}
\DeclareMathSymbol{\epsilon}{\mathalpha}{letters}{"0F}
\DeclareMathSymbol{\zeta}{\mathalpha}{letters}{"10}
\DeclareMathSymbol{\eta}{\mathalpha}{letters}{"11}
\DeclareMathSymbol{\theta}{\mathalpha}{letters}{"12}
\DeclareMathSymbol{\iota}{\mathalpha}{letters}{"13}
\DeclareMathSymbol{\kappa}{\mathalpha}{letters}{"14}
\DeclareMathSymbol{\lambda}{\mathalpha}{letters}{"15}
\DeclareMathSymbol{\mu}{\mathalpha}{letters}{"16}
\DeclareMathSymbol{\nu}{\mathalpha}{letters}{"17}
\DeclareMathSymbol{\xi}{\mathalpha}{letters}{"18}
\DeclareMathSymbol{\pi}{\mathalpha}{letters}{"19}
\DeclareMathSymbol{\rho}{\mathalpha}{letters}{"1A}
\DeclareMathSymbol{\sigma}{\mathalpha}{letters}{"1B}
\DeclareMathSymbol{\tau}{\mathalpha}{letters}{"1C}
\DeclareMathSymbol{\upsilon}{\mathalpha}{letters}{"1D}
\DeclareMathSymbol{\phi}{\mathalpha}{letters}{"1E}
\DeclareMathSymbol{\chi}{\mathalpha}{letters}{"1F}
\DeclareMathSymbol{\psi}{\mathalpha}{letters}{"20}
\DeclareMathSymbol{\omega}{\mathalpha}{letters}{"21}
\DeclareMathSymbol{\varepsilon}{\mathalpha}{letters}{"22}
\DeclareMathSymbol{\vartheta}{\mathalpha}{letters}{"23}
\DeclareMathSymbol{\varpi}{\mathalpha}{letters}{"24}
\DeclareMathSymbol{\varrho}{\mathalpha}{letters}{"25}
\DeclareMathSymbol{\varsigma}{\mathalpha}{letters}{"26}
\DeclareMathSymbol{\varphi}{\mathalpha}{letters}{"27}

\RedeclareMathOperator{\div}{\textbf{div}}

\renewcommand{\epsilon}{\varepsilon}

\newcommand{\q}{\bm q}

\renewcommand{\r}{\vec r}
\newcommand{\f}{\bm f}
\newcommand{\g}{\bm g}

\newcommand{\F}{\vec F}
\newcommand{\G}{\vec G}

\newcommand{\df}[2]{\frac{\partial #1}{\partial #2}}
\newcommand{\con}{\q}
\newcommand{\prim}{\vec u}

\newcommand{\dt}{\Delta t}

\newcommand{\fl}{\f}
\newcommand{\gl}{\g}

\newcommand{\half}{\frac{1}{2}}

\renewcommand{\vec}[1]{\mathbf{#1}}

\renewcommand{\div}{\ensuremath{\textbf{div}}\xspace}
\renewcommand{\epsilon}{\varepsilon}

\newcommand{\order}{\ensuremath{\mathcal{O}}}
\newcommand{\oo}[1]{\ensuremath{\frac{1}{#1}}}

\newcommand{\x}{\ensuremath{\vec x}\xspace}

\newcommand{\vel}{\ensuremath{\vec v}\xspace}

\newcommand\eg{e.g.\ }

\renewcommand{\dt}{\partial_t}

\newcommand{\dx}{\partial_x}
\newcommand{\dy}{\partial_y}

\newcommand{\evar}{\vec r}

\newcommand{\M}{\mathcal M}

\newcommand{\new}[1]{#1\xspace}

\let\origdoublepage\cleardoublepage
\renewcommand{\cleardoublepage}{%
	\clearpage{\pagestyle{empty}\origdoublepage}}

\makeatletter
\newcommand{\ChapterOutsidePart}{%
	\def\toclevel@chapter{-1}\def\toclevel@section{0}\def\toclevel@subsection{1}}
\newcommand{\ChapterInsidePart}{%
	\def\toclevel@chapter{0}\def\toclevel@section{1}\def\toclevel@subsection{2}}
\makeatother

\begin{document}


\begin{frontmatter}

\title{Entropy Stable Numerical Fluxes for Compressible Euler Equations which are Suitable for All Mach Numbers
	}
\author[1]{Jonas P.\ Berberich\corref{cor1}}
\ead{jonas.berberich@mathematik.uni-wuerzburg.de}
\author[1]{Christian Klingenberg}
\cortext[cor1]{Corresponding author: 
	Tel.: +49 931 31-88861;  
	Fax: +49 931 31-83494;}
\ead{klingen@mathematik.uni-wuerzburg.de}
\address[1]{Dept.~of Mathematics, Univ.~of W\"urzburg, Emil-Fischer-Stra{\ss}e 40, 97074 W\"urzburg, Germany}

\begin{abstract}
	We propose two novel two-state approximate Riemann solvers for the compressible Euler equations which are provably entropy dissipative and suitable for the simulation of low Mach numbers. What is new, is that one of our two methods in addition is provably kinetic energy stable. Both methods are based on the entropy satisfying and kinetic energy consistent methods of Chandrashekar (2013). The low Mach number compliance is achieved by rescaling some speed of sound terms in the diffusion matrix in the spirit of Li \& Gu (2008). In numerical tests we demonstrate the low Mach number compliance and the entropy stability of the proposed fluxes.
\end{abstract}


\end{frontmatter}



\section{Introduction}
Compressible Euler equations are used to model the flow of compressible inviscid fluids such as air. To find approximate solutions of the Euler system it is common to use finite volume methods. These are well-suited due to their conservative nature and their capability to resolve discontinuities. The fluxes at the cell interfaces are often determined by approximating the solution of the 1-d interface Riemann problem using numerical (two-state) fluxes.

For a sequence of ever lower Mach numbers, the solutions of the compressible Euler equations with well-prepared initial data converge towards solutions of the incompressible Euler equations \cite{Feireisl2019}. This limit, however, is not correctly represented in a finite volume scheme using conventional numerical fluxes due to excessive diffusion at low Mach numbers. Special low Mach number compliant numerical fluxes have been developed (e.g.\ \cite{Turkel1999,Li2008,Dellacherie2010,Rieper2011,Miczek2015,Oswald2015,Chalons2016,Barsukow2017,Zenk2017}, to correct this behavior.

Stability is required to ensure the convergence of a finite volume method. There are different notions of stability, one of them being entropy stability. Entropy stability is a non-linear stability criterion which additionally ensures that the entropy inequality -an admissibility criterion for physical solutions- is satisfied. Ismail \& Roe \cite{Ismail2009}, and later Chandrashekar \cite{Chandrashekar2013}, developed two-state Riemann solvers based on the Roe flux \cite{Roe1981}, which ensure entropy dissipation to achieve entropy stability. The flux by Chandrashekar \cite{Chandrashekar2013} is kinetic energy consistent additionally.

Recently, a numerical flux based on \cite{Ismail2009} has been developed which is entropy dissipative and low Mach compliant \cite{Chen2018}. In this article we present two methods \new{for compressible Euler equations closed with an ideal gas law} based on the entropy stable and kinetic energy compliant fluxes of Chandrashekar \cite{Chandrashekar2013} and a low Mach modification of Li \& Gu \cite{Li2008}. The methods we propose are entropy dissipative and low Mach compliant. One of our methods is additionally kinetic energy stable. \new{The low Mach number method from Li \& Gu \cite{Li2008} has the tendency to develop so-called checkerboard instabilities at low Mach numbers. Numerical experiments indicate that the methods developed in this article suppress the checkerboard modes.}

The rest of the article is structured as follows:
In \cref{sec:EulerEquations} we describe the entropy and kinetic energy dissipative methods introduced in \cite{Chandrashekar2013}. In \cref{sec:theoLowmach} we discuss the behavior of the method at low Mach numbers and we present a correction in \cref{sec:lowMach_fix}. Numerical tests demonstrating the improved results at low Mach numbers and the entropy dissipation of the proposed methods are presented in \cref{sec:tests}.

\section{Kinetic Energy and Entropy Stable Fluxes}
\label{sec:EulerEquations}
The 2-d Euler equations which model the conservation laws of mass, momentum, and energy of a compressible inviscid fluid are given by
\begin{equation}
\df{\con}{t} + \df{\fl}{x} + \df{\gl}{y} = 0,
\label{eq:eul2d_cartesian}
\end{equation}
where the conserved variables and fluxes are
\begin{equation}
\con = \begin{bmatrix}
\rho \\ \rho u \\ \rho v \\ E \end{bmatrix}, \quad
\fl = \begin{bmatrix}
\rho u \\
\rho u^2+p \\
\rho u v \\
(E+p)u \end{bmatrix}, \quad 
\gl = \begin{bmatrix}
\rho v \\
\rho u v \\
\rho v^2+p \\
(E+p)v \end{bmatrix}.
\end{equation}
Moreover, $E=\rho \varepsilon + \tfrac{1}{2}\rho|\vec v|^2 $ is the total energy per unit volume with $\vec v=[u,v]^T$ being the velocity. The pressure $p$ is related to the density and internal energy via the ideal gas equation of state
\begin{equation}
p=R T\rho\qquad\text{with}\qquad T=\frac{\gamma-1}{R}\frac\epsilon\rho.
\end{equation}
The dependent variable $T$ is called temperature, the constants are the gas constant $R$, and the ratio of specific heats $\gamma$.

\subsection{Entropy-Entropy Flux and Entropy Variables}
A pair $(U,\bm\phi)$ with a convex function $U(\con)$ and a vector valued function $\bm\phi(\con)=[\phi_x(\con),\phi_y(\con)]^T$ is called \emph{entropy-entropy flux pair}, if it satisfies the relations
\begin{equation}
U^\prime(\con)\f^\prime(\con) = \phi_x^\prime(\con),\qquad
U^\prime(\con)\g^\prime(\con) = \phi_y^\prime(\con).
\end{equation}
Using this pair, we can add the additional conservation law 
\begin{equation}
\df{U}{t} + \df{\phi_x}{x}+ \df{\phi_y}{y}= 0
\end{equation}
to \cref{eq:eul2d_cartesian}. As usual in the context of hyperbolic conservation laws, we also want to admit discontinuous solutions and interpret all the derivatives in \cref{eq:eul2d_cartesian} in the weak sense. At discontinuities, the entropy is not necessarily conserved. Instead, the inequality
\begin{equation}
\df{U}{t} + \df{\phi_x}{x}+ \df{\phi_y}{y}\leq 0
\label{eq:entropy_condition}
\end{equation}
is demanded as a criterion to choose admissible (physical) solutions.
We define \emph{entropy variables} by
\begin{equation}
\evar(\con):=U^\prime(\con)
\end{equation}
and the dual to the entropy flux by $\bm\psi(\evar)=[\psi_x(\evar),\psi_y(\evar)]^T$ with
\begin{equation}
\psi_x(\evar):=\evar\cdot\f(\con(\evar)) - \phi_x(\con(\evar)),\quad
\psi_y(\evar):=\evar\cdot\g(\con(\evar)) - \phi_y(\con(\evar)),
\end{equation}
where $\con(\evar) $ is the inverse of $\evar(\con)$ defined above. The inverse exists because of the convexity of $U(\con)$.

For the Euler equations \eqref{eq:eul2d_cartesian}, the most common choice of an entropy-entropy flux pair is
\begin{equation}
U:=-\frac{\rho s}{\gamma-1}, \qquad \bm\phi:=-\frac{\rho \vel s}{\gamma-1},
\end{equation}
where $s:=\ln\left(p\rho^{-\gamma}\right) = -(\gamma-1)\ln(\rho)-\ln(\beta) -\ln(2)$ up to a constant with $\beta:=1/(2RT)$. This choice is not unique \cite{Harten1983d}, but it is the one consistent with the entropy condition from thermodynamics in the presence of heat transfer \cite{Hughes1986}. The entropy variables are subsequently given by
\begin{equation}
\evar:=\begin{bmatrix}
\frac{[\gamma-s]}{\gamma-1}-\beta |\vel|^2,&
2\beta u,&
2\beta v,&
-2\beta
\end{bmatrix}^T.
\end{equation}
and the entropy flux dual
\begin{equation}
\bm\psi=\rho \vel.
\end{equation}

\subsection{A Basic Finite Volume Method}
In this section, for brevity reasons, we only describe a simple quadrature-free finite volume method, which is a second order accurate finite volume method on a static Cartesian grid. In practice, more elaborate methods are used (see e.g.\ \cite{Toro2009}).

We divide the domain $\Omega=[a,b]\times[c,d]$ with $a<b,c<d$ into cells 
\begin{equation}
\Omega_{ij}:=\left[x_{i-\half} ,x_{i+\half}\right]\times\left[y_{j-\half} ,y_{j+\half}\right]
\end{equation}
for $i=0,\dots,N-1$, $j=0,\dots,M-1$. The cell-interface centers are 
\begin{equation}
\x_{i-\half,j}:=\left[a+i\Delta x,c+\left(j+\half\right)\Delta y\right]^T, \x_{i,j-\half}:=\left[a+\left(i+\half\right)\Delta x,c+j\Delta y\right]^T
\end{equation}
with $\Delta x:=\frac{b-a}N, \Delta y:=\frac{d-c}M$.
We integrate $\con$ in each cell to obtain cell-averaged values
\begin{equation}
\hat\con_{ij}(t):=\frac1{\Delta x\Delta y}\int_{\Omega_{ij}}  \con(\x,t) \,d\x.
\end{equation}
We find an evolution equation for the cell-average values by cell-wise integrating \cref{eq:eul2d_cartesian}:
\begin{equation}
\dt\hat\con_{ij}(t) =- \frac1{\Delta x\Delta y}\int_{\Omega_{ij}}\dx\f(\con(\x,t))+\dy\g(\con(\x,t))\,d\x.
\end{equation}
To construct our simple finite volume method we use Fubini's theorem, the fundamental theorem of calculus, and an approximation of the interface integral by the interface centered point value. The interface centered fluxes are approximated using a numerical two-state flux. This yields
\begin{align}
\dt\hat\con_{ij}(t) \approx &- \frac1{\Delta x}\left[\F\left(\hat\con^-_{i+\half,j}(t),\hat\con^+_{i+\half,j}(t)\right)-\F\left(\hat\con^-_{i-\half,j}(t),\hat\con^+_{i-\half,j}(t)\right)\right]\nonumber\\
&- \frac1{\Delta y}\left[\G\left(\hat\con^-_{i,j+\half}(t),\hat\con^+_{i,j+\half}(t)\right)-\G\left(\hat\con^-_{i,j-\half}(t),\hat\con^+_{i,j-\half}(t)\right)\right],
\label{eq:semidiscrete}
\end{align}
where the $\hat{\con}^\pm$ values are obtained using a non-oscillatory reconstruction on the cell-average values.
This set of ODEs is then integrated numerically to evolve the approximate solution $\hat\con$ in time.
In the rest of the article we drop the hat at $\hat{\con}_{ij}$ and just write $\con_{ij}$. Also, in the rest of the article we will only consider numerical fluxes in $x$-direction, since the fluxes $\f$ and $\g$ for Euler equations can be converted into each other by only correctly rotating velocity vectors. For symmetry reasons we assume $\F$ and $\G$ to have the same relation.

\subsection{Entropy Conservative Numerical Fluxes}
Tadmor \cite{Tadmor1987,Tadmor2003} introduced the concept of \emph{entropy conservative numerical fluxes} $\F^\text{ec}$, which have to satisfy the relation
\begin{equation}
(\evar(\con^+)-\evar(\con^-))\cdot\F^\text{ec}(\con^-,\con^+)=\psi(\evar(\con^+))-\psi(\evar(\con^-))=(\rho\vel)^+-(\rho\vel)^-.
\end{equation}
The last identity is only valid for the Euler equations \eqref{eq:eul2d_cartesian}.
Different entropy conservative numerical fluxes have been proposed by Tadmor \cite{Tadmor1987}, Ismail \& Roe \cite{Ismail2009}, and Chandrashekar \cite{Chandrashekar2013}. Our method is based on the numerical flux by Chandrashekar \cite{Chandrashekar2013}, which can be written as
\begin{equation}
\F^\ast(\con^-,\con^+):=
\begin{bmatrix}
F^{\ast,\rho}\\
F^{\ast,\rho u}\\
F^{\ast,\rho v}\\
F^{\ast,E}
\end{bmatrix}:=
\begin{bmatrix}
\hat{\rho}\bar u\\
\bar u F^{\ast,\rho}+\tilde p\\
\bar v F^{\ast,\rho}\\
\left(\frac1{2(\gamma-1)\hat\beta}-\frac12\overline{|\vel|^2}\right)F^{\ast,\rho}
+\bar{\vel}\cdot[F^{\ast,\rho u},F^{\ast,\rho v}]^T
\end{bmatrix}.
\label{eq:Praveen_central_1d}
\end{equation}
The averages are the arithmetic average
$\bar a:=\frac12(a^-+a^+)$ and the logarithmic average $\hat a := \frac{a^+-a^-}{\ln a^+-\ln a^-}$. A non-singular implementation of $\hat a$ is presented in \cite{Ismail2009}. The pressure average is $\tilde p:=\bar{\rho}/(2\bar{\beta})$ where $\bar\beta$ is computed from $\beta^\pm=\rho^\pm/(2 p^\pm)$. This pressure average corresponds to the harmonic average in the temperature \cite{Chandrashekar2013}. The notations for the averages are used throughout this article. 

\subsection{Kinetic Energy Preserving Fluxes}
From the density and momentum equations in the Euler equations \eqref{eq:eul2d_cartesian} the balance law 
\begin{equation}
\df K t + \df {(K u)}x + \df {(K v)}y = -u\df p x  -v\df p y
\end{equation}
for the kinetic energy $K:=\frac12\rho |\vel|^2$ can be derived. Integration over the whole domain $\Omega$ while ignoring the boundaries yields
\begin{equation}
\df{}{t}\int_{\Omega} K\,d\x= \int_{\Omega}p\df u x +p \df v y \,d\x.
\label{eq:total_kinetic_energy_dissipation}
\end{equation}
Jameson \cite{Jameson2008} shows that any numerical flux $\F^J$ which can be formulated in the form
\begin{equation}
\F^\text{J}\left(\con^-,\con^+\right)=
\begin{bmatrix}
F^{\text{J},\rho}\\
\bar u F^{\text{J},\rho} + \left<p\right>\\
\bar v F^{\text{J},\rho}\\
F^{\text{J},E}
\end{bmatrix},
\label{eq:J_flux}
\end{equation}
satisfies the discrete analogon 
\begin{align}
\df{}{t}\sum_{i,j}K_{ij}\Delta x \Delta y  
&= \sum_{i,j}\left(-\half |\vel_{ij}|^2\df{\rho_{ij}}{t}+\vel_{ij}\df{(\rho\vel)_{ij}}{t}\right) \Delta x\Delta y\nonumber\\
&= \sum_{i,j} \left(\left<p\right>_{i+\half,j}\frac{\Delta u_{i+\half,j}}{\Delta x}+\left<p\right>_{i,j+\half}\frac{\Delta v_{i,j+\half}}{\Delta y}\right) \Delta x\Delta y
\label{eq:discrete_total_kinetic_energy_dissipation}
\end{align}
of \cref{eq:total_kinetic_energy_dissipation}. \Cref{eq:discrete_total_kinetic_energy_dissipation} can easily be computed using the flux from \cref{eq:J_flux} and the corresponding flux in $y$-direction $\G^J$. The fluxes $F^{\text{J},\rho}$ and $F^{\text{J},E}$ are consistent approximations of the density and energy flux and $\left< p\right>$ approximates the interface pressure. Clearly, the entropy conservative numerical flux $\F^\ast$ from \cref{eq:Praveen_central_1d} is in this kinetic energy preserving form.

\new{In the low Mach number limit, the right-hand side of \cref{eq:total_kinetic_energy_dissipation} vanishes with the divergence of velocity (\eg \cite{Guillard1999}). \Cref{eq:total_kinetic_energy_dissipation} then describes the conservation of kinetic energy. This makes kinetic energy consistency especially relevant for low Mach number fluxes. Most numerical flux functions violate this condition. For example, in \cite{Miczek2015} it is numerically shown that the kinetic energy rises for a simulation of the incompressible Gresho \cite{Gresho1990} vortex using a central flux and implicit time stepping.}

\subsection{Entropy Diffusion}
\label{sec:entropy_diffusion}
For the scheme to be stable in the presence of discontinuities it needs to dissipate entropy. Following \cite{Chandrashekar2013} this is achieved by modifying the Roe scheme diffusion \cite{Roe1981} such that the Roe matrix is applied to the jump in entropy variables $\evar$
instead of conserved variables $\con$. The standard Roe scheme uses the diffusion
\begin{equation}
\F^\text{Roe,diff}\left(\con^-,\con^+\right):= -\half D^\text{Roe}\Delta\con:= -\half R|\Lambda|^\text{Roe}R^{-1}\Delta\con,
\end{equation}
where 
\begin{equation}
R:=\begin{bmatrix}
1              & 1             & 0         & 1             \\
u-c            & u             & 0         & u+c           \\
v              & v             & -1        & v             \\
H-cu           & \half|\vel|^2 &  -v       & H+cu
\end{bmatrix},
\end{equation}
with the enthalpy $H=\frac{c^2}{\gamma-1}+\frac{|\vel|^2}2$, is the matrix of right eigenvectors of $\df{\f(\con)}{\con}$ and
\begin{equation}
|\Lambda|^\text{Roe}:=
\text{diag}\left(\begin{bmatrix}
|\lambda_1|,
|\lambda_2|,
|\lambda_3|,
|\lambda_4|
\end{bmatrix}\right)=\text{diag}\left(\begin{bmatrix}
|u-c|,
|u|,
|u|,
|u+c|
\end{bmatrix}\right)
\label{eq:matrix_of_eigenvalues}
\end{equation}
is the diagonal matrix with the absolute values of the corresponding eigenvalues. The whole matrix $D^\text{Roe}$ is evaluated at the Roe average state \cite{Roe1981} to ensure accurate shock capturing. In order to apply the diffusion matrix $D^\text{Roe}$ to the jump in entropy variables, we have to transform these to conserved variables. So the diffusion part of our numerical flux is
\begin{equation}
\F^\text{ES,diff}\left(\con^-,\con^+\right) := -\half R|\Lambda|^\text{Roe}R^{-1}\frac{\partial \con}{\partial \r}\Delta\r.
\end{equation}
The entropy diffusion $\F^\text{ES,diff}$ can be formulated in a simpler form which can lead to a more efficient implementation: It is shown by Barth \cite{Barth1999} that there is a scaling $\tilde R=R S^{-\half}$ of the Eigenvectors $R$ with $\frac{\partial \con}{\partial \r}=\tilde R \tilde R^T$ which leads to the form
\begin{equation}
\F^\text{ES,diff}\left(\con^-,\con^+\right)=-\half \tilde R|\Lambda|^\text{Roe}\tilde R^{-1}\tilde R \tilde R^T\Delta \r = -\half \underbrace{R|\Lambda|^\text{Roe}SR^T}\limits_{=:Q}\Delta \r
= -\half Q\Delta \r
\end{equation}
with the scaling matrix
\begin{equation}
S:=\text{diag}\left(\begin{bmatrix}
\frac{\rho}{2\gamma},&
\frac{(\gamma-1)\rho}{\gamma},&
p,&
\frac{\rho}{2\gamma}
\end{bmatrix}\right).
\end{equation}
Since $Q$ is positive definite by construction, $\F^\text{ES,diff}$ is dissipative in the entropy variables.
The numerical flux
\begin{equation}
\F^\text{ES}\left(\con^-,\con^+\right) := \F^\ast\left(\con^-,\con^+\right) + \F^\text{ES,diff}\left(\con^-,\con^+\right)
\end{equation}
is hence entropy satisfying\new{\ in the sense, that a spatially discrete analogon of the entropy inequality is satisfied. Proofs are analogously to \cite{Chandrashekar2013,Ray2013}}. The numerical flux in $y$-direction is constructed in the same way.

\subsection{Kinetic Energy Diffusion}
Alongside entropy stability we also aim for kinetic energy stability. We have seen that $\F^\ast$ ist consistent with the evolution of kinetic energy derived from the Euler equations \eqref{eq:eul2d_cartesian}. In order to guarantee kinetic energy stability, the diffusive term in the numerical flux has to do dissipate kinetic energy. Chandrashekar \cite{Chandrashekar2013} shows that this requires the condition $|\lambda_1|=|\lambda_4|$ to hold in \cref{eq:matrix_of_eigenvalues}. The most obvious way to achieve this is to choose
\begin{equation}
|\Lambda|^\text{KES}:=
\text{diag}\left(\begin{bmatrix}
|\lambda_1|,
|\lambda_2|,
|\lambda_3|,
|\lambda_4|
\end{bmatrix}\right)=
\text{diag}\left(\begin{bmatrix}
|u|+c,
|u|,
|u|,
|u|+c
\end{bmatrix}\right).
\label{eq:matrix_of_eigenvalues_KES}
\end{equation}
in the definition of the entropy-diffusion matrix $Q$. The numerical flux defined by this modification on the ES scheme will be called ES-KES or $\F^\text{ES-KES}$ throughout this article. Note that this method is more diffusive than the ES scheme. However, it adds one more relevant stability property.

\subsection{Intermediate State}
\label{sec:intermediate_state}
To ensure correct upwinding, the diffusion matrix in the standard Roe scheme is evaluated at the so-called Roe average state. Anyway, in our numerical flux $\F^\text{ES}$ we use $\F^\ast$ instead of the standard central flux, so we can not expect the shock-capturing property to still hold for our method. In the following we discuss at which intermediate state the diffusion matrix $Q$ should be evaluated.

The entropy stability property does not depend on the intermediate state, since $Q$ is positive definite for any state by construction.
The kinetic energy stability also does not depend on the particular choice of the intermediate state, only on the relation of the entries in the diagonal matrix $|\Lambda|$.
For the flux to have the contact property, Chandrashekar \cite{Chandrashekar2013} shows that we have to choose
\begin{equation}
c_\text{int}=\sqrt{\frac{\gamma}{2\hat{\beta}}} \quad\text{and}\quad H_\text{int}=\frac{c_\text{int}^2}{\gamma-1}+\frac12 |\vel_\text{int}|^2.
\label{c_average}
\end{equation}
All other averages can be chosen freely, so we can use arithmetic or logarithmic averages for example.

However, for implementation reasons it can be useful to hand a intermediate state vector in primitive variables to the routine which computes the diffusion matrix. We can realize this using the average state
\begin{align}
\vel_\text{int}&:=\bar \vel&
p_\text{int}&:=\bar p&
\rho_\text{int}&:=2 p_\text{int} \hat \beta.
\end{align}
In the computation of the diffusion matrix we compute all the other variables from the primitive intermediate state in the straight forward way, e.g. $c_\text{int}=\sqrt{\gamma p_\text{int}/\rho_\text{int}}$ and the enthalpy as described in \cref{c_average}.

\section{Low Mach Number Asymptotics}
\label{sec:theoLowmach}
The Euler equations \eqref{eq:eul2d_cartesian} can be cast in the non-dimensional form
\begin{equation}
\df{}{t}\begin{bmatrix}
\rho \\ \rho u \\ \rho v \\ E \end{bmatrix} 
+ \df{}{x}\begin{bmatrix}
\rho u \\
\rho u^2 + \oo{\M^2}p \\
\rho u v \\
(E+p)u \end{bmatrix} 
+ \df{}{y}\begin{bmatrix}
\rho v \\
\rho u v \\
\rho v^2 + \oo{\M^2}p \\
(E+p)v \end{bmatrix} = 0
\label{eq:eul2d_nondim}
\end{equation}
using only the assumption that the reference velocity is computed as the quotient of the reference length and time and one parameter, which we will call reference Mach number $\M$.
%

The low Mach number limit of the Euler equations \cref{eq:eul2d_nondim} is well-known and studied (\eg \cite{Guillard1999,Barsukow2017b,Feireisl2019}). For well-prepared initial data a series of solutions of \cref{eq:eul2d_nondim} with different reference Mach numbers converge to solutions of the incompressible Euler equations for $\M\to 0$. 
Conventional finite volume methods tend to fail to correctly represent this limit for their numerical solutions. One reason for this is excessive diffusion at low Mach numbers. 

Consider the ES-flux introduced in \cref{sec:entropy_diffusion}. Note that, again, we only consider the flux in $x$-direction for simplicity. For small jumps we can approximate
\begin{align}
\F^\text{ES}\left(\con^-,\con^+\right) 
&= \F^\ast\left(\con^-,\con^+\right) - \frac12D^\text{Roe}(\con_\text{int})\left.\df{\con(\evar)}{\evar}\right|_{\evar=\evar\left(\con_\text{int}\right)}\left(\evar\left(\con^+\right)-\evar\left(\con^-\right)\right)\nonumber\\
&\approx J\left(\con_\text{int}\right)\con_\text{int} - \frac12D^\text{Roe}\left(\con_\text{int}\right)\Delta\con\nonumber\\
&= \left.\df{\con}{\prim}\right|_{\con_\text{int}}J_\text{prim}\left(\con_\text{int}\right)
\left.\df{\prim}{\con}\right|_{\con_\text{int}}
\con_\text{int} - \frac12\left.\df{\con}{\prim}\right|_{\con_\text{int}}D^\text{Roe}_\text{prim}\left(\con_\text{int}\right)\left.\df{\prim}{\con}\right|_{\con_\text{int}}\Delta\con
\label{eq:flux_approximation}
\end{align}
because of $\con^-\approx\con_\text{int}\approx\con^+$ and consequently $\F^\ast\left(\con^-,\con^+\right)\approx\f(\con_\text{int})$. The intermediate state $\con_\text{int}$ is the one defined in \cref{sec:intermediate_state}. The flux Jacobian in primitive variables is $J_\text{prim}=\df{\prim}{\con}J\df{\con}{\prim}$ with the flux jacobian in conserved variables $J=\df{\f}{\con}$ and the Roe diffusion matrix in primitive variables is $D^\text{Roe}_\text{prim}=\df{\prim}{\con}D^\text{Roe}\df{\con}{\prim}$. The primitive variables are
\begin{equation}
\prim:=[\rho,u,v,p]^T.
\end{equation}

\Cref{eq:flux_approximation} justifies the comparison of the two matrices $J_\text{prim}$ and $D^\text{Roe}_\text{prim}$ with regard of their formal scaling with the reference Mach number to gain insight into the asymptotic behavior of the method for small Mach numbers:

\begin{align}
J_\text{prim}&=
\left[\begin{array}{*{8}c}
\mathcal{O}\left(1\right)      & \mathcal{O}\left(1\right)      & 0                              & 0                              \\
0                              & \mathcal{O}\left(1\right)      & 0                              & \mathcal{O}\left(\frac{1}{\M^{2}}\right) \\
0                              & 0                              & \mathcal{O}\left(1\right)      & 0                              \\
0                              & \mathcal{O}\left(1\right)      & 0                              & \mathcal{O}\left(1\right)      \\
\end{array}\right],\label{eq:scalings_jacobian}\\
D^\text{Roe}_\text{prim}&=\left[\begin{array}{*{8}c}
\mathcal{O}\left(1\right)      & 0     & 0 & \mathcal{O}\left(\frac{1}{\M}\right) \\
0                              & \mathcal{O}\left(\frac{1}{\M}\right) & 0     & \mathcal{O}\left(\frac{1}{\M}\right) \\
0                              & 0     & \mathcal{O}\left(1\right)      & \mathcal{O}\left(1\right)      \\
0                              & 0     & 0 & \mathcal{O}\left(\frac{1}{\M}\right) \\
\end{array}\right]+\order\left(\M\right).
\label{eq:scalings_roe}
\end{align}
Note that the diffusion matrix scaling for the ES-KES flux is the same, since for low Mach numbers $|u+c|\approx|u|+c\approx|u-c|$.
From \cref{eq:scalings_jacobian,eq:scalings_roe} we see that there are some terms in which the diffusion matrix dominates the flux Jacobian for small Mach numbers. This explains the excessive diffusion of the method at low Mach numbers. Different methods have been proposed to correct those terms in the Roe diffusion matrix (e.g. \cite{Turkel1999,Rieper2011,Miczek2015,Oswald2015,Barsukow2017}). For the entropy stability of our method, however, the modification should keep the positive definiteness of $Q$. The following modification provides this.

\section{Low Mach Modifications of the ES and ES-KES Fluxes}
\label{sec:lowMach_fix}

Following \cite{Li2008} we modify the diagonal matrices \cref{eq:matrix_of_eigenvalues,eq:matrix_of_eigenvalues_KES} to make the schemes low Mach number compliant: We use 
\begin{align}
|\Lambda|^\text{Roe}_\text{LM}&:=\text{diag}\left(\begin{bmatrix}
|u-\tilde c|,
|u|,
|u|,
|u+\tilde c|
\end{bmatrix}\right),\label{eq:MLmod_roe}\\
|\Lambda|^\text{KES}_\text{LM}&:=\text{diag}\left(\begin{bmatrix}
|u|+\tilde c,
|u|,
|u|,
|u|+\tilde c
\end{bmatrix}\right)
\label{eq:MLmod_kes}
\end{align}
with 
\begin{equation}
\tilde c := c\cdot \max(\min(M,1),M_\text{cut})\quad\text{with}\quad M:=\frac{|\vel|}{c}\quad \text{and}\quad M_\text{cut}\in[0,1].
\end{equation}
Other possible definitions of the rescaled speed of sound $\tilde c$ can be found e.g.\ in \cite{Li2008,Chen2018}.
The global cut-off Mach number $M_\text{cut}$ can be used to increase the diffusion and thus the stability at Mach numbers which are lower than the one expected in a particular simulation.
Using those diagonal matrices instead of $|\Lambda|^\text{Roe}$ and $|\Lambda|^\text{KES}$ yields the  numerical fluxes $\F^\text{ES-LM}$ and $\F^\text{ES-KES-LM}$. It is obvious that this modification is compatible with the proof of entropy stability in \cite{Chandrashekar2013}. Also, $\F^\text{ES-KES-LM}$ still satisfies the relation $|\lambda_1|=|\lambda_4|$ required for kinetic energy stability.
This modification changes \cref{eq:scalings_roe} into
\begin{equation}
\left(D^\text{Roe}_\text{LM}\right)_\text{prim}=\left[\begin{array}{*{8}c}
\mathcal{O}\left(1\right)      & 0    & 0 & \mathcal{O}\left(1\right)      \\
0                              & \mathcal{O}\left(1\right)      & 0 & \mathcal{O}\left(\frac{1}{\M}\right) \\
0                              & 0 & \mathcal{O}\left(1\right)      & 0     \\
0                              & 0     & 0 &\mathcal{O}\left(1\right)      \\
\end{array}\right]+\order\left(\M\right).
\label{eq:scalings_roe_lm}
\end{equation}
and the terms of dominating diffusion at low Mach numbers vanish. As before, this scaling is also valid for the diffusion matrix of the ES-KES-LM flux.

\section{Numerical Tests}
\label{sec:tests}
In our tests we include the standard Roe flux \cite{Roe1981} (Roe), the Roe flux with the low Mach fix \cite{Li2008} described in \cref{eq:MLmod_roe} (Roe-LM), the entropy stable fluxes (ES, ES-KES) from \cite{Chandrashekar2013} and the entropy stable fluxes with low Mach fix (ES-LM, ES-KES-LM) proposed in this article. In the low Mach methods we use $M_\text{cut}=0$. For time-stepping we use the third order accurate four stage Runge--Kutta Method as described in \cite{Kraaijevanger1991}. In practice, an implicit method should be used to evolve low Mach number flows in time due to the stiffness in time \cite{Miczek2015}. However, this is not in the scope of this short article.

For the equation of state we choose $\gamma=1.4$, which is a suitable value to describe air. We set the value of the gas constant to $R=1$, such that density and pressure have the same order of magnitude.

\subsection{Test of the Entropy Stability}
\label{sec:test_entropy}
\begin{figure}[t]
	\centering
	\includegraphics[width=\textwidth]{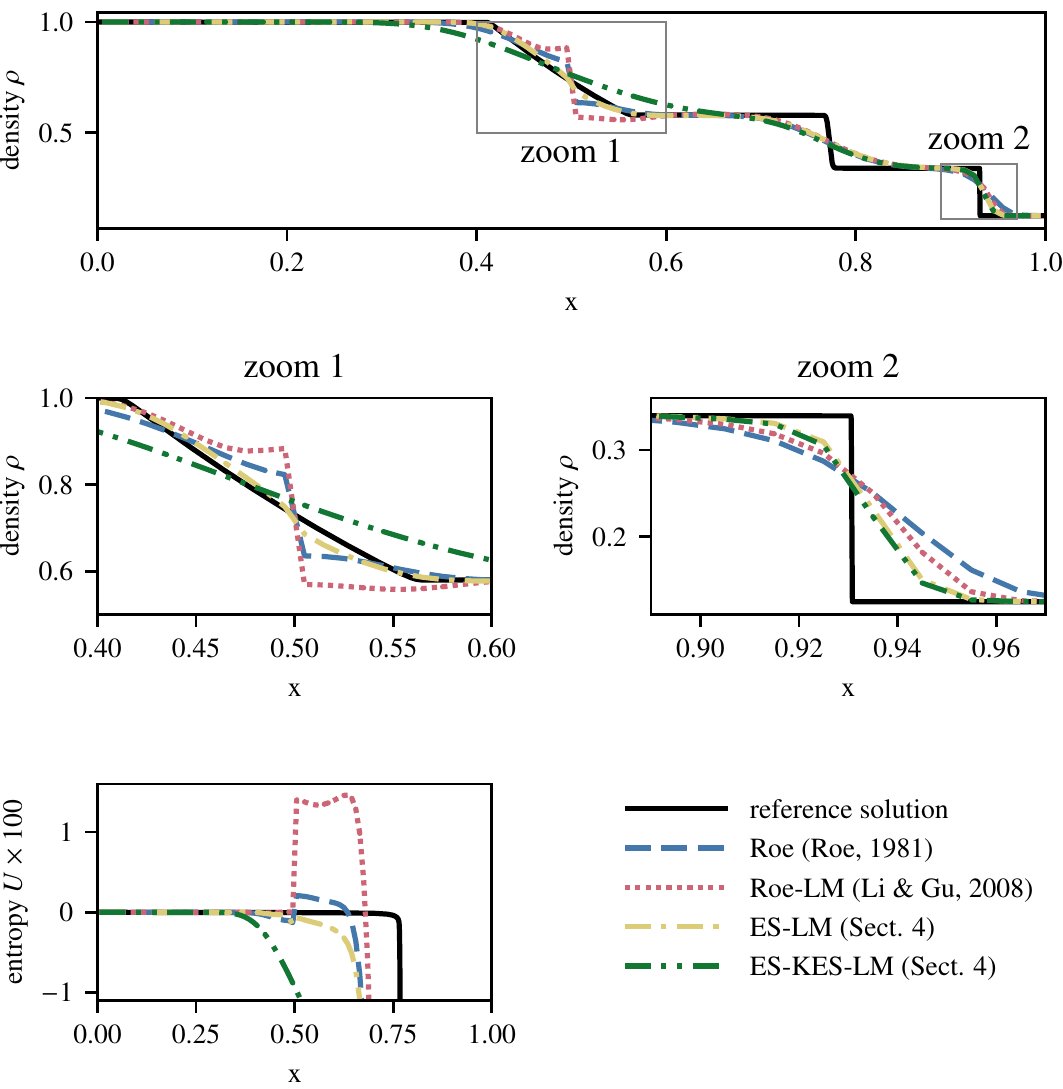}
	\caption{\label{fig:entropy_Sod}Standing sound wave test at $t=0.2$ with different numerical fluxes at a 100 cells grid. The test-setup is described in \cref{sec:test_entropy}. Top: density on the whole domain. Middle: Magnified views of particular regions of the top panel. Bottom: Entropy at the whole domain}
\end{figure}

One test for the entropy compliance of a method is a standing sound wave, for which e.g. the Roe scheme is well-known to produce a non-physical jump. We solve the 1-d Riemann-problem given by
\begin{align}
(\rho, u, p)(x<0.5)&:=(1,0.75,1)\\
(\rho, u, p)(x\geq 0.5)&:=(0.125,0,0.1)	
\end{align}
on 100 equidistant grid cells in the domain $\Omega=[0,1]$. We use the standard Roe-flux, the standard Roe-flux with the low Mach modification \cref{eq:MLmod_roe} (Roe-LM), and the entropy stable low Mach methods ES-LM and ES-KES-LM introduced in this article (\cref{sec:lowMach_fix}) with constant (which means no) reconstruction. The result at $t=0.2$ is shown in \cref{fig:entropy_Sod}.
As a reference solution we use a simulation with the local Lax--Friedrichs numerical flux on $100\,000$ grid cells. We see the non-entropic jump which is produced by the Roe scheme. The low Mach modification in Roe-LM even increases the non-entropic jump. As expected, for the entropy stable all Mach methods ES-LM and ES-KES-LM there is no significant non-entropic jump. It is notable that the methods proposed in this article also have an improved accuracy on the expansion shock (zoom 2). In the bottom panel of \cref{fig:entropy_Sod} we see the entropy of the solution at time $t=0.2$. While the entropy at initial time is non-positive on the whole domain (This is easy to compute from the initial conditions), the Roe and Roe-LM method lead to positive values of entropy at time $t=0.2$. The ES-LM and ES-KES-LM method lead to non-positive entropy values on the whole domain.

\subsection{Test of the Contact Property}
\label{sec:test_contact}
\begin{figure}[t]
	\centering
	\includegraphics[width=\textwidth]{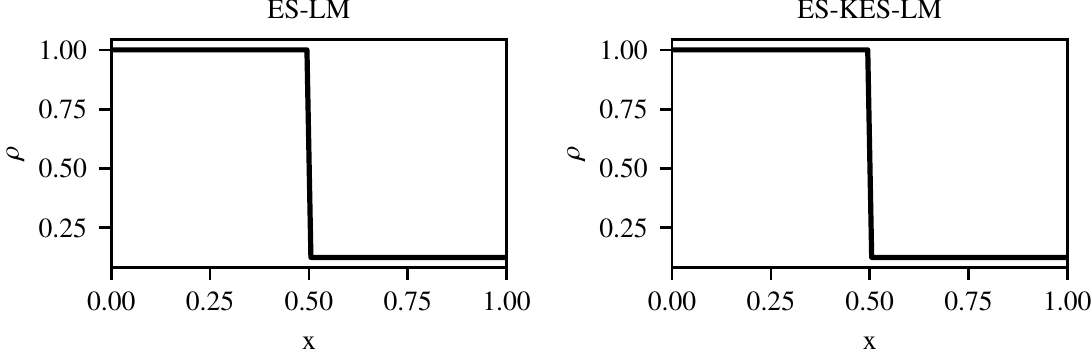}
	\caption{\label{fig:contact}Density of a contact discontinuity at $t=0.2$ computed using the entropy stable low Mach methods ES-LM and ES-KES-LM. The setup is described in \cref{sec:test_contact}}
\end{figure}

We test the contact property of the method using the Riemann problem
\begin{align}
(\rho, u, p)(x<0.5)&:=(1,0.75,1)\\
(\rho, u, p)(x\geq 0.5)&:=(0.125,0,0.1)	
\end{align}
on 100 equidistant grid cells in the domain $[0,1]$.  We use the ES-LM and ES-KES-LM fluxes with constant reconstruction and the intermediate state described in \cref{sec:intermediate_state}. The result at $t=0.2$ is shown in \cref{fig:contact}. Both methods accurately resolve the contact discontinuity.

\subsection{Low Mach Gresho Vortex}
\label{sec:gresho}

\begin{figure}
	\centering
	\includegraphics[scale=1]{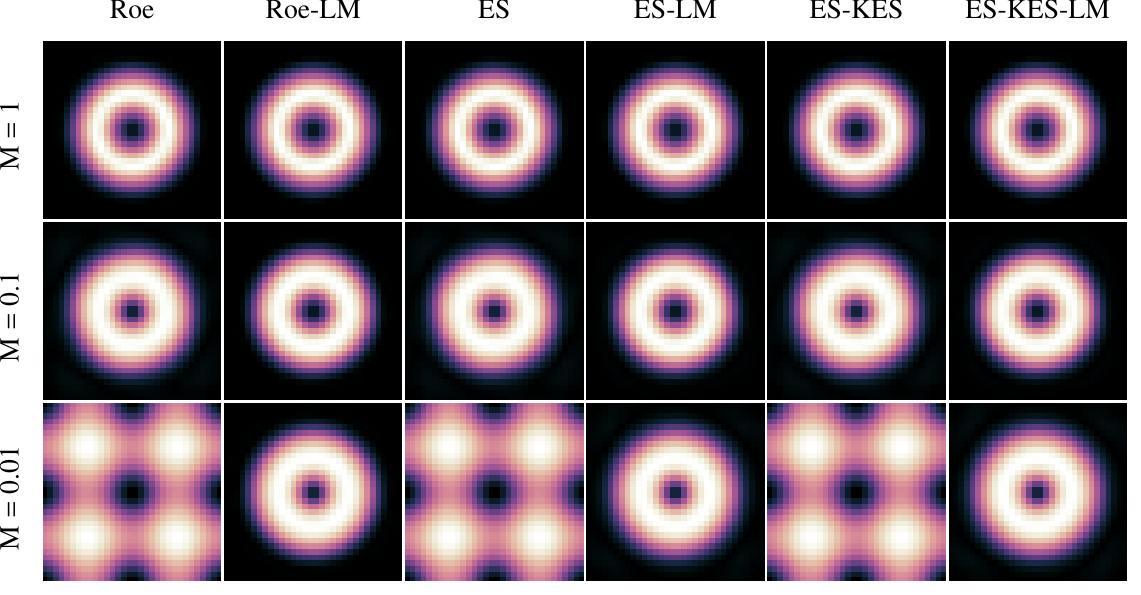}
	\caption{\label{fig:gresho}Local Mach number of the Gresho vortex test from \cref{sec:gresho} after 0.1 rotations. The initial maximal Mach number decreases from top to bottom. In the different columns different numerical fluxes are applied}
\end{figure}

The incompressible Gresho vortex \cite{Gresho1990} can be extended to a family of stationary solutions of the Euler equations with a parameter that adjusts the maximal local Mach number in the setup \cite{Barsukow2017}. The setup is
\begin{equation*}
(\rho,\vel,p) = 
\begin{cases}
(1,\,5 r \vec e_\phi,\,p_c + \frac{25}2 r^2)&\text{if } r<0.2,\\
\left(1,\,(2-5 r)\vec e_\phi, \,p_c+4\ln(5r)+4-20 r +\frac{25}2 r^2\right) &\text{if } 
0.2\leq r<0.4,\\
(1,\,0,\,p_c +4\ln(2)-2) & \text{else},
\end{cases}
\end{equation*}
where $p_c =\frac{1}{\gamma \tilde M^2}
\frac12$ and $\vec e_\phi$ is the unit vector in angular direction.
We apply the standard Roe flux, the Roe flux with low Mach fix \cite{Li2008}, the entropy stable fluxes from \cite{Chandrashekar2013} and the entropy stable fluxes with low Mach fix proposed in this article. We use limited linear reconstruction on a $32\times32$ cells grid with periodic boundary conditions to evolve the vortex for 0.1 revolutions. For the low Mach number fluxes we use $M_\text{cut}=0$. The Mach number at final time is shown in \cref{fig:gresho} for the different numerical fluxes (columns) and maximal initial Mach number parameters $\tilde M=1,0.1,0.01$ (rows). The Roe, ES, and ES-KES numerical fluxes lead to completely diffused vortexes for lower Mach number, while the proposed fluxes (ES-LM, ES-KES-LM) are as capable of accurately resolving the Gresho vortex at lower Mach numbers as the Roe method with low Mach fix (Roe-LM, \cite{Li2008}).


\section{Conclusions and Outlook}
\label{sec:conclusions}
We presented novel numerical flux functions which combine the entropy and kinetic energy stability properties of the fluxes proposed in \cite{Chandrashekar2013} with the low Mach number compliance of the method from \cite{Li2008}. The entropy stability and the low Mach number compliance have been shown in numerical tests. The contact property holds due to the correct choice of the intermediate state. It is worth noting that the proposed methods also show an improved performance at fast shocks. For practical applications of the proposed methods we suggest the combination with implicit time-stepping to overcome the stiffness in time. \new{To extend the method to higher order while maintaining the entropy stability, special care has to be given to the reconstruction procedure, as discussed in \cite{Fjordholm2012,Ray2016}. To avoid the carbuncle phenomenon at strong shocks, \cite{Chandrashekar2013,Ray2013} suggest a hybrid diffusion with a Rusanov-type diffusion term.}

\section*{Acknowledgments}
The authors thankfully acknowledge the helpful discussions with Praveen Chandrashekar.
Jonas Berberich thanks the organizers of the NumHyp 2019 for a science-wise interesting and socially pleasurable conference.
The research of Jonas Berberich is supported by the Klaus Tschira Foundation.

\bibliographystyle{model1-num-names}
\bibliography{library}
\end{document}